\documentclass[12pt,twoside]{amsart}

\usepackage{graphicx} 
 \long\def\red#1{{\color{red}#1}}

\usepackage{amsmath,amsthm,amssymb,amscd}
\usepackage{enumerate}
\usepackage{latexsym}
\usepackage[top=1.5in, bottom=1.5in, left=1.25in, right=1.25in]	 {geometry}

\usepackage{tikz}
\usetikzlibrary{arrows, decorations.pathmorphing, backgrounds, positioning, fit, petri}
\usetikzlibrary{ decorations.text}

\def\f{\frac}

\def\({\left(}
\def \){ \right)}
\def\[{\left[}
\def \]{ \right]}
\def\Bl{\Bigl}
\def\Br{\Bigr}

\def\sa{{\sigma}}

 \def\s{{\sigma}}
 
 \def\va{\varepsilon}

  \def\sph{\mathbb{S}^{d}}

\def\ld{\lambda}
\def\bl{\bigl}
\def\br{\bigr}

\newcommand{\wh}{\widehat}

\usepackage{tikz}\usepackage{tikz-3dplot}

\usepackage{mathtools}
\mathtoolsset{showonlyrefs,showmanualtags}

\usepackage{hyperref} 
\hypersetup{
    colorlinks=true,       
    linkcolor=blue,          
    citecolor=magenta,        
    filecolor=magenta,      
    urlcolor=cyan           
}

\theoremstyle{plain} 
\newtheorem{lemma}[equation]{Lemma}
\newtheorem{proposition}[equation]{Proposition}
\newtheorem{theorem}[equation]{Theorem}
\newtheorem{corollary}[equation]{Corollary}

\theoremstyle{definition}

\theoremstyle{remark}
\newtheorem{remark}[equation]{Remark}

\numberwithin{equation}{section}

\begin{document}

\title{Stolarsky  principle and energy optimization on the sphere}

\author{Dmitriy Bilyk}
\address{School of Mathematics, University of Minnesota, Minneapolis, MN 55408, USA.}
\email{dbilyk@math.umn.edu}

\author{Feng Dai}
\address{Department of Mathematical and Statistical Sciences\\
University of Alberta\\ Edmonton, Alberta T6G 2G1, Canada.}
\email{fdai@ualberta.ca}

\author{Ryan Matzke}
\address{School of Mathematics, University of Minnesota, Minneapolis, MN 55408, USA.}
\email{matzk053@umn.edu}

\thanks{This work is partially supported by the Simons foundation collaboration grant (Bilyk),   NSERC  Canada under
grant  RGPIN 04702 (Dai), and the NSF Graduate Research  Fellowship (Matzke). The first two authors are grateful to CRM Barcelona: this collaboration  originated  during their participation in the research program ``Constructive Approximation and Harmonic Analysis (Bilyk's trip was sponsored by NSF grant DMS 1613790). }

\maketitle

\begin{abstract}
The classical {\it{Stolarsky invariance principle}} connects the spherical cap $L^2$ discrepancy of a finite point set on the sphere to the pairwise sum of Euclidean distances between the points. In this paper we further explore and extend this phenomenon. In addition to a new elementary proof of this fact, we establish several new analogs, which relate various notions of discrepancy to different discrete energies. In particular, we find that the hemisphere discrepancy is related to the sum of geodesic distances. We also extend these results to arbitrary measures on the sphere and arbitrary notions of discrepancy and apply them to problems of  energy optimization and combinatorial geometry and find that, surprisingly, the geodesic distance energy behaves differently than its Euclidean counterpart.  
\end{abstract}

\section{Introduction and main results}

In numerous areas of mathematics and other sciences, one is faced with the problem of distributing large finite sets of points on the sphere as uniformly as possible. There exist various quantitative measures of uniformity of spherical point distributions. Among the most popular ones are discrepancy and energy.

Let $Z = \{ z_1,..., z_N\} \subset \sph$ be an $N$-point set in the $d$-dimensional sphere, and let $\sigma$ denote the normalized Lebesgue surface measure on $\mathbb S^d$, i.e. $\sigma (\sph) =1$. For a given subset of the sphere, $A\subset \sph$, the discrepancy of $Z$ with respect to $A$ is defined as
\begin{equation}\label{1disc}
 D( Z, A) = \frac{1}{N}  \sum_{k=1}^N {{\bf 1}}_A (z_k)- \sigma (A),
 \end{equation}
  in other words, $D(Z,A)$ indicates how well the Lebesgue measure of $A$ is approximated by the counting measure $\frac1{N} { \sum_{k=1}^N } \delta_{z_k}$. To obtain good  finite distributions $Z$, one normally evaluates and strives to minimize the supremum (extremal discrepancy) or average (e.g., $L^2$ discrepancy) of $|D(Z,A)|$ over some rich and well-structured collection of sets $A$. Typical examples of such collections include spherical caps, slices, convex sets etc. -- specific choice depends on the problem at hand. For a good exposition of discrepancy and numerical integrations on the sphere, as well as  discrepancy theory in general, the reader is referred to, e.g., \cite{DT,M99}.

On the other hand, the energy of $Z$ with respect to a function $F: \sph \times \sph \rightarrow {\mathbb{R}}$  is defined as
\begin{equation}\label{1energy} E_F (Z) = \frac1{N^2} \sum_{i,j = 1}^N F (z_i, z_j),
\end{equation} i.e. points of $Z$ are viewed as interacting particles on the sphere which repel according to the potential given by $F$. If $F (x,x)$ is undefined, which is the case for the most common example arising in electrostatics,   the Riesz potential $F(x,y) = \| x - y \|^{-s}$, the diagonal terms are omitted in the sum above. In many situations, minimizing (or maximizing, depending on the structure of $F$) the energy $E_F (Z)$  yields well-distributed  point-sets on the sphere, and the quality of this distribution may be measured by the difference of the discrete energy $E_F (Z)$ and the energy of the continuous uniform distribution $I_F (\sigma) = \int_{\sph} \int_{\sph}  F(x,y) d\s (x) d\s (y)$. Vast literature exists on problems of this nature: we refer the reader, e.g., to   the upcoming book \cite{BHSbook}.

It is known that in some cases these two ways of quantifying  equidistribution are closely connected. One of the first instances of such a connection was obtained in 1973 by Stolarsky \cite{St}, who proved that minimizing the $L^2$ discrepancy with respect to spherical caps is equivalent to maximizing the pairwise sum of Euclidean distances, i.e. $E_F (Z)$ with $F(x,y) = \| x-y \|$. More precisely, he established the identity
$$ c_d  [ D_{L^2,\textup{cap}} (Z) ]^2   =   \int\limits_{\mathbb S^{d}} \int\limits_{\mathbb S^{d}} \|  x- y \| \, d\sigma (x)\, d\sigma (y)\,\,  - \,\, \frac{1}{N^2} \sum_{i,j = 1}^N \| z_i - z_j \|,   $$ which came to be known as the {\emph{Stolarsky invariance principle}} (see \S \ref{ClSt} for more details).\\ 

In this paper we further explore these connections and obtain various  versions of the Stolarsky principle in different settings, yielding new relations between discrepancy and energy optimization.

In \S \ref{ClSt} we revisit the classical Stolarsky invariance principle  and give a very simple elementary proof of this identity.

In \S \ref{HemSt} we observe that replacing {\emph{all}} spherical caps with  {\emph{hemispheres}} (i.e. spherical caps with aperture $\pi/2$) one obtains a variant of Stolarsky principle with the geodesic distance $d(x,y)$ in place of  the Euclidean distance, Theorem \ref{SPhem}. This allows one to easily  characterize finite point sets on $\sph$, which maximize the sum of geodesic distances in all dimensions $d\ge 1$ (for  even $N$ these are just symmetric sets), and points to a drastic difference with the case of  Euclidean distances, see Theorem \ref{t.geodpoints}. 

In \S \ref{GeodInt} we take this idea one step further and show that an analog of the  Stolarsky principle holds for general probability measures $\mu$ in place of the counting measure (Theorem \ref{SPmu}). This provides a way to characterize the maximizers of the geodesic distance energy integral $\displaystyle{ I_g (\mu ) = \int_{\sph} \int_{\sph} \big( d(x,y) \big)^\delta d\mu (x) d\mu (y)}$, with $\delta >0$,  over all probability measures on the sphere, Theorem \ref{GeodMax}. In particular,
\begin{itemize}
\item [$\circ$] for $0< \delta <1$, the unique maximizer is $\sigma$;
\item [$\circ$] for $\delta =1$, maximizers are centrally symmetric measures;
\item [$\circ$] for $\delta >1$, maximizers are measures of the form $\mu = \frac12 (\delta_{p} + \delta_{-p})$.
\end{itemize}
The second part of this statement is a  consequence  of Stolarsky principle. The case $\delta \in (0,1)$  is actually proved in the   companion paper of the authors \cite{BD2}, which also studies the cases of $\delta < 0$ (geodesic Riesz energy) and $\delta = 0$ (logarithmic energy) by means of  analyzing ultraspherical expansions. This brings up a surprising difference between the geodesic and Euclidean settings. In the latter case a result of  Bjorck \cite{Bjorck} from 1955, see Theorem \ref{Bjo}, states that the critical value is $\delta =2$ rather than $\delta =1$. This effect in dimension $d=1$ has been previously noticed in \cite{BHS}.

In \S \ref{GenSt} we explore the connections between energy optimization and discrepancy on a more general level. We show that for positive definite functions $F$ one can define a natural notion of discrepancy, so that the analog of the Stolarsky principle holds for general measures $\mu$, Theorem \ref{thm-3-2}, and make connections between positive definiteness and refinements of this property and  minimization of energy integrals $I _F (\mu) = \int_{\sph} \int_{\sph} F(x\cdot y) d\mu (x) d\mu (y)$, in particular, the question  whether $\sigma$ is a unique  minimizer of $I_F (\mu)$, in other words, whether attaining  equilibrium under the potential $F$ imposes uniform distribution. These results are further developed and applied   in the parallel paper of the authors \cite{BD2}, yielding sharp asymptotic estimates for the difference  $E_F (Z) - I_F (\sigma)$ as well as new proofs of classical discrepancy bounds  \eqref{Beck}.

In \S \ref{WedSli}, we consider   discrepancy with respect to spherical ``slices" (i.e. intersections of half-spaces; previously studied in \cite{Blu}) and spherical ``wedges" (i.e. symmetrized slices -- this notion of discrepancy came up recently \cite{BLd} in connection to problems of uniform tessellation of the sphere by hyperplanes and one-bit compressed sensing).  In these cases, one also obtains Stolarsky-type identities with potentials $F(x,y) = \big( 1 - d(x,y) \big)^2$ and $F(x,y) = \big( \frac12 - d(x,y) \big)^2$, respectively.

Finally, in the appendix, \S \ref{App}, we compute the values of some of the spherical  integrals that arise in the exposition.\\

We note that the basic strategy behind  most  versions of Stolarsky principle, at a very low level, is straightforward. Computing the $L^2$ discrepancy, one  squares out the expression in \eqref{1disc}, thus pairwise interactions between points of $Z$ arise from cross terms of the form ${\bf 1}_A (z_i) \cdot {\bf 1}_A (z_j)$. When integrated over the the test sets $A$ in a given class, this  yields the interaction potential $F(z_i, z_j)$, which is often represented as the volume of intersection of test sets ``centered" at $z_i$ and $z_j$,  see e.g. \eqref{L2ta}, however, the details in some settings get rather technical. This approach is employed in \S \ref{ClSt}-\ref{GeodInt},  \S \ref{WedSli}. A similar idea has been used by Torquato \cite{Torq} for ``number variance", a quantity very similar to $L^2$ discrepancy.   In \S \ref{GenSt} we go in the opposite direction and show that for any positive definite interaction potential one can construct an appropriate notion of discrepancy, so that the Stolarsky principle holds.

We would also like to mention that the interest in Stolarsky principle in different settings has recently spiked: \cite{BD}    studied it from the point of view of numerical integration on the sphere, \cite{Owen} uses it in  applications to genomics, \cite{Skr} explores Stolarsky principle in general metric spaces, \cite{BLd}  connects it to tessellations of the sphere, while the present paper and \cite{BD2} deal with it in the context of energy optimization.

In the text,  the dimension of the sphere is $d \ge 1$, i.e. $\mathbb S^d \subset \mathbb R^{d+1}$; $\sigma$ is the {\emph{normalized}}   Lebesgue ($d$-dimensional Hausdorff) measure on $\mathbb S^d$, i.e. $\sigma (\sph) =1$; $\|x \|$ denotes the Euclidean norm of $x\in \mathbb R^{d+1}$; and $d(x,y) $ denotes the   geodesic distance between $x$ and $y \in \mathbb S^d$, {\emph{normalized}} so that the distance between antipodal poles is equal to $1$, i.e. $d(x,y) = \frac1{\pi} \cos^{-1} ( x \cdot y )$. The cardinality of a finite set $Z$ is denoted by $\# Z$.  The set of all finite signed Borel measures on $\sph$ is denoted by  $\mathcal B$, and $\mathcal M$ stands for the set of  probability measures on $\sph$ (positive Borel measures with total mass one). Further explanations, background information and references will be given in each section.

\section{Classical Stolarsky invariance principle for spherical caps.}\label{ClSt}

We consider  ``spherical caps" $C(x,t)$ with center $x\in \mathbb S^d$ and ``height" $t
\in [-1, 1]$, i.e.
\begin{equation}
C(x,t ) = \{ z \in \mathbb S^d:\, z \cdot x > t \} .
\end{equation}
We define the $L^2$ discrepancy of $Z$ with respect to spherical caps:
\begin{equation}
[ D_{L^2,\textup{cap}} (Z) ]^2  =  \int\limits_{-1}^1  \int\limits_{\mathbb S^{d}}  \bigg| \frac1{N} \sum_{j=1}^N {\bf 1}_{C(x,t)} (z_j) -  \sigma \big( C(x,t) \big) \bigg|^2 d\sigma (x) \, dt = \int\limits_{-1}^1  [ D^{(t)}_{L^2,\textup{cap}} (Z) ]^2  \, dt .
\end{equation}
The  following result was proved by Stolarsky in 1973 \cite{St}: 

\begin{theorem}[Stolarsky invariance principle] \label{t.SP}
Let $Z = \{ z_1 , z_2, \dots, z_N \} \subset \mathbb S^d$.
Then the following relation holds
\begin{equation}\label{SP}
   [ D_{L^2,\textup{cap}} (Z) ]^2   = C_d\,  \Bigg(\, \, \int\limits_{\mathbb S^{d}} \int\limits_{\mathbb S^{d}} \|  x- y \| \, d\sigma (x)\, d\sigma (y)\,\,  - \,\, \frac{1}{N^2} \sum_{i,j = 1}^N \| z_i - z_j \| \Bigg).
\end{equation}
The constant $C_d $ satisfies
\begin{equation}
C_d = \frac12 \int_{\mathbb S^{d}} |p\cdot z| \, d\sigma (z) = \frac1{d} \frac{\omega_{d-1}}{\omega_{d}} = \frac{v_{d}}{\omega_{d}} = \frac{1}{d} \frac{\Gamma \big((d+1)/2\big)}{\sqrt{\pi} \, \Gamma(d/2)} \sim \frac{1}{\sqrt{2\pi d} }\,\,\, \textup{ as } \,\, d\rightarrow \infty,
\end{equation}
where $\omega_{d}$ is the surface area of $\mathbb S^{d}$ and $v_{d}$ is the volume of the unit ball in $\mathbb R^{d}$.\\

\end{theorem}

\noindent This theorem states that
\begin{itemize}
\item minimizing the $L^2$ spherical cap discrepancy of $Z$  is equivalent to maximizing the sum of Euclidean distances between the points of $Z$.

\item the $L^2$ spherical cap discrepancy can be realized as the difference between the continuous and discrete energies or, equivalently, the error of numerical integration of the distance  integral $\int\limits_{\mathbb S^{d}} \int\limits_{\mathbb S^{d}} \|  x- y \| \, d\sigma (x)\, d\sigma (y)$ by the cubature  formula with knots at the points of $Z$.
\end{itemize}

It is well known \cite{Beck1,Beck2}  that the optimal order of the $L^2 $ spherical cap discrepancy is $N^{-\frac12 - \frac1{2d}}$, i.e.
\begin{equation}\label{Beck}
 c_d N^{-\frac12 - \frac1{2d}} \le \inf_{
 \# Z = N }  D_{L^2,\textup{cap}} (Z)   \le c'_d N^{-\frac12 - \frac1{2d}},
\end{equation}
which in turn  bounds  the difference of continuous and discrete energies in \eqref{SP}.


In addition to the original proof in \cite{St}, an alternative proof has been given in \cite{BD}. 
Here we present  a new short and simple proof of the Stolarsky invariance principle \eqref{SP}. It strongly resonates with the proof in \cite{BD}, but is completely elementary in nature. A similar  proof in a probabilistic interpretation has been independently given in \cite{Owen} (compare Lemmas \ref{Lm1} and \ref{Lm2} below to Proposition 1 of \cite{Owen}), and analogous ideas are used in \cite{Skr} on general metric spaces.

The proof  of \eqref{SP} follows the aforementioned strategy: one  squares out the integrand, and the discrete part (pairwise interactions) arises naturally from the cross terms. The important ingredient is the  following relation between intersections of spherical caps and the Euclidean distance between their centers: 
\begin{lemma}\label{Lm1} For arbitrary  $x$, $y\in \mathbb S^{d}$ we have
\begin{equation}\label{inter}
 \int\limits_{-1}^1 \sigma \big( C(x,t) \cap C(y,t) \big) \, dt  = 1- C_{d} \| x- y\|,
\end{equation}
where the constant $C_d$ is given by $C_d = \frac12 \int\limits_{\mathbb S^{d}} |p\cdot z| \, d\sigma (z)$ for any fixed point $p\in \mathbb S^{d}$.
\end{lemma}

\begin{proof} Recall that $\sigma$ is normalized so that  $\sigma ({\mathbb S^{d} })  = 1$. We have
\begin{align*}
 \int\limits_{-1}^1 \sigma \big( C(x,t)  \cap C(y,t) \big) \, dt &  = \int\limits_{-1}^1 \int\limits_{\mathbb S^{d}} \mathbf 1_{C(x,t)} (z) \cdot  \mathbf 1_{C(y,t)} (z) d\sigma (z) \, dt \\
&  =  \int\limits_{\mathbb S^{d}} \int\limits_{-1}^1 \mathbf 1_{C(z,t)} (x) \cdot  \mathbf 1_{C(z,t)} (y) \, dt \, d\sigma (z) \,  =  \int\limits_{\mathbb S^{d}} \int\limits_{-1}^{ \min \{ x\cdot z, y\cdot z  \}}  \, dt \, d\sigma (z) \, \\
& = \int\limits_{\mathbb S^{d}}  \big(  \min \{ x\cdot z, y\cdot z  \}  + 1 \big) \, d\sigma (z)\red{.}
\end{align*}
We now write $ \min \{ x\cdot z, y\cdot z  \}  = \frac12 \big(   x\cdot z + y\cdot z  - \big| (x-y) \cdot z \big|\big)$. Obviously  $\int\limits_{\mathbb S^{d}} x\cdot z \,d\sigma(z)  =  \int\limits_{\mathbb S^{d}}  y\cdot z \,d\sigma (z) = 0$. By rotational invariance, we observe that $$  \int\limits_{\mathbb S^{d}} \big| (x-y) \cdot z \big| \, d\sigma (z) =  \| x- y\| \cdot \int\limits_{\mathbb S^{d}} \left| \frac{ x-y }{\| x - y\|} \cdot z \right| \, d\sigma (z) = 2C_d \| x- y \|,$$
where $C_d = \frac12 \int\limits_{\mathbb S^{d}} |p\cdot z| \, d\sigma (z)$, and this finishes   the proof.
\end{proof}

Here we  essentially repeated the proof from \cite{BD}, but the proof of the next lemma, which gives the quadratic mean value of the size of the spherical caps, is simpler (does not use reproducing kernels).
\begin{lemma}\label{Lm2}
For any $p \in \mathbb S^{d}$ we have
\begin{equation}\label{quadave}
\int\limits_{-1}^1 \left( \sigma \big( C (p,t) \big) \right)^2 \, dt = 1 - C_d \int\limits_{\mathbb S^{d}} \| x - p \| \,d\sigma (x),
\end{equation}
\end{lemma}

\noindent {\bf{Remark:}} It is clear that $\displaystyle{\int\limits_{\mathbb S^{d}} \| x - p \|\,  d\sigma (x) = \int\limits_{\mathbb S^{d}} \int\limits_{\mathbb S^{d}} \| x - y \| \,d\sigma (x)\, d\sigma (y).}$

\begin{proof}
We use the result of the previous lemma, i.e. relation \eqref{inter}, to compute
\begin{align*}
1 - C_d\int\limits_{\mathbb S^{d}} \| x - p \| \,d\sigma (x) & = \int\limits_{\mathbb S^{d}}\big(  1-  C_d \| x - p \| \big) \, d\sigma (x) \\
& = \int\limits_{\mathbb S^{d}} \int\limits_{-1}^1 \int\limits_{\mathbb S^{d}} \mathbf 1_{C(x,t)} (z) \cdot  \mathbf 1_{C(p,t)} (z) d\sigma (z) \,  d t\,  d\sigma (x) \\
&=  \int\limits_{-1}^1 \int\limits_{\mathbb S^{d}} \mathbf 1_{C(p,t)} (z)  \left( \,\,\, \int\limits_{\mathbb S^{d}} \mathbf 1_{C(z,t)} (x)  \, d\sigma (x) \right)\,  d\sigma (z)\, dt \\
&  = \int\limits_{-1}^1 \left( \sigma \big( C (p,t) \big) \right)^2 \, dt .
\end{align*}
\end{proof}


\noindent This is one of numerous examples of a situation in which averaging over scales  simplifies things. For the $L^2$ discrepancy for spherical caps of fixed height $t$:
\begin{equation}\label{L2t}
D^{(t)}_{L^2,\textup{cap}} (Z) := \left(   \int_{\mathbb S^{d}}  \bigg| \frac1{N}  \sum_{j=1}^N {\bf 1}_{C(x,t)} (z_j) -  \sigma \big( C(x,t) \big) \bigg|^2 d\sigma (x) \right)^{1/2},
\end{equation}
one would have to deal  with   $\sigma \big( C(x,t) \cap C(y,t) \big)$, which  has complicated structure,  and no short   relation akin to \eqref{inter} is available, see e.g. \cite{Owen}. Hence in this case, there is no formula as succinct and explicit  as  the Stolarsky principle, however one can still write  down a generic relation where the interactions between $z_j$'s would involve $ \sigma \big( C(z_i,t) \cap C(z_j,t) \big)$. 
\begin{proposition} For any $Z =\{z_1,z_2, \dots, z_N \} \subset \mathbb S^d$ and a fixed $t\in [-1,1]$,  the following relation holds
\begin{equation}\label{L2ta}
\left[ D^{(t)}_{L^2, \textup{cap}} (Z) \right]^2  = \frac1{N^2} \sum_{i,j=1}^N \sigma  \big( C(z_i,t) \cap C(z_j,t) \big) - \left( \sigma \big( C(p,t)\big) \right)^2 ,
\end{equation}
where $p$ is an arbitrary point on the sphere.
\end{proposition}




\begin{proof} We note that $\sigma \big(C (x,t) \big)$ is independent of $x\in \mathbb S^d$, hence
\begin{align}\label{L2tb}
\left[ D^{(t)}_{L^2} (Z) \right]^2  & =  \int\limits_{\mathbb S^{d}} \bigg|  \frac1{N} \sum_{j=1}^N {\bf 1}_{C(x,t)} (z_j) -  \sigma \big( C(x,t) \big) \bigg|^2 d\sigma (x)  \\
\nonumber & = \frac1{N^2}  \sum_{i,j=1}^N \int\limits_{\mathbb S^{d}} \mathbf 1_{C(x,t)} (z_i) \cdot  \mathbf 1_{C(x,t)} (z_j) d\sigma (x) \\
\nonumber &  \phantom{ = } \,\,  - \frac{2}{N}  \int\limits_{\mathbb S^{d}} \sum_{j=1}^N {\mathbf 1}_{C(z_j,t)} (x) \cdot \sigma \big( C(x,t)\big) \, d\sigma (x) \\
\nonumber & \phantom{ = } \,\, + \left( \sigma \big( C(p,t)\big) \right)^2 \\
\nonumber &=  \frac1{N^2} \sum_{i,j=1}^N \sigma  \big( C({z_i},t) \cap C({z_j},t) \big) - \left( \sigma \big( C(p,t)\big) \right)^2 .
\end{align}
\end{proof}

 Integrating   identity \eqref{L2ta} with respect to $t$ and applying relations \eqref{inter} and \eqref{quadave} we finish the proof of  the Stolarsky principle \eqref{SP}.


\section{Stolarsky principle for hemispheres.}\label{HemSt}

An (open) hemisphere in the direction of $x\in \mathbb S^d$ is simply a spherical cap of height $t=0$:
\begin{equation}
H( x) = \{ z \in \mathbb S^d:\, z\cdot x >0 \}  = C(x,0).
\end{equation}
Since $\sigma \big(H(x)\big) = \frac12$, the  natural $L^2$ discrepancy for this set system  is
\begin{equation}\label{e.dischem}
D_{L^2,\textup{hem}} (Z) :=  \left(   \int_{\mathbb S^{d}}  \bigg| \frac1{N}  \sum_{j=1}^N {\bf 1}_{H(x)} (z_j) -  \frac12 \bigg|^2 d\sigma (x) \right)^{1/2} =  D^{(0)}_{L^2,\textup{cap}} (Z)\red{.}
\end{equation}
While, as mentioned above, generally the quantity $\sigma \big( C(x,t) \cap C(y,t) \big)$ is complicated, in the case $t=0$ (hemispheres) it has a very simple representation: for $x$, $y \in \mathbb S^d$
\begin{equation}\label{2hem}
\sigma \big(H(x) \cap H(y) \big)  =  \sigma \big( C(x,0) \cap C(y,0) \big) = \frac12 \cdot \Big( 1 - d(x,y)  \Big),
\end{equation}
where $d(x,y)$ is the normalized geodesic distance on the sphere between $x$ and $y$. This can be very easily seen from the picture below.

\begin{center}
\begin{tikzpicture}[scale=0.75]
\draw (0,0) circle [radius = 4cm];
\draw[dashed] ( 3.46cm, -2cm) -- (-5.19 cm ,3 cm) ;
\draw[dashed] ( -4cm, 0) -- (6cm,0) ;
\filldraw[fill=cyan, draw=black] (0,0) -- ( 4cm, 0) arc (0: 150: 4cm) -- (0,0) ;
\filldraw (0, 4cm) circle (2pt) node[above] {$x$};
\filldraw (2cm, 3.46cm) circle (2pt) node[above] {$y$};
\def\myshift#1{ \raisebox{ 1ex}}
\draw[ <->, postaction= {decorate, decoration= {text along path, text align= center, text= {|\sffamily\myshift|{$d(x,y)$}}}} ] (0, 4.7cm) arc (90: 60: 4.7cm);
\draw[ <->, postaction= {decorate, decoration= {text along path, text align= center, text= {|\sffamily\myshift|{ $\sigma   = \frac{1}{2} (1 -d(x,y)) $}}}} ] (-5.19 cm ,3 cm) arc (150: 0: 6 cm);
\end{tikzpicture}
\end{center}


\noindent Combining \eqref{2hem}, Proposition \ref{L2t}  and the obvious fact that    $ \int\limits_{\mathbb S^{d}} \int\limits_{\mathbb S^{d}} d( x, y)  \, d\sigma (x)\, d\sigma (y)  = \frac12$,  
one  arrives at the following result.

\begin{theorem}[Stolarsky principle for hemispheres]\label{SPhem} For any $N$-point set $Z \in \mathbb S^d$, the following relation holds:
\begin{equation}\label{e.SPhem}
[D_{L^2,\textup{hem}} (Z)]^2      = \frac12 \left(\, \int\limits_{\mathbb S^{d}} \int\limits_{\mathbb S^{d}} d( x, y)  \, d\sigma (x)\, d\sigma (y) - \frac{1}{N^2} \sum_{i,j = 1}^N d(z_i , z_j) \right).
\end{equation}
\end{theorem}

The statement looks strikingly similar to the original Stolarsky principle \eqref{SP}. One can say that the Euclidean distance corresponds to the \emph{mean} over $t\in [-1,1]$, while the geodesic distance corresponds to the \emph{median} ($t=0$) of the heights of the spherical caps. Despite the fact that the original Stolarsky principle has been proved in 1973, this version is new and  has not been observed before. (At the time of preparation of this manuscript we have learned that this version of the Stolarsky principle has been independently and simultaneously proved by Skriganov \cite{Skr}.)  Relation \eqref{e.SPhem} has several interesting features and consequences.


First of all, the principle of \emph{irregularities of distribution} does not hold in this situation, that is, the hemisphere discrepancy can be very small, even zero, for large $N$. Indeed,  for any symmetric distribution $Z$, it is easy to see that the $L^2$ hemisphere discrepancy is equal to zero.  Moreover, \eqref{e.SPhem} allows us to characterize finite point distributions in $\mathbb S^d$, which maximize the sum of geodesic distances. 
\begin{theorem}\label{t.geodpoints} Let $d\ge 1$. Then the following holds:
\begin{enumerate}[(i)]
\item\label{M1} For any point distribution $Z = \{ z_1, \dots, z_N \} \subset \mathbb S^d$,
\begin{equation}
\frac{1}{N^2} \sum_{i,j = 1}^N d(z_i , z_j)  \le \frac12.
\end{equation}
\item\label{M2} For a given $N \in \mathbb N$ the sum above is maximized if and only if  the following condition holds: for any $x\in \mathbb S^d$, such that the hyperplane  $x^\perp$ contains no points of $Z$, the numbers of points of $Z$ on either side of $x^\perp$ differ by at most one, i.e.
 \begin{equation}
\big| \# \big(  Z \cap H(x) \big)  - \# \big( Z \cap H(-x) \big) \big| \le 1.
\end{equation}
\item\label{M3} If $N$  is even,  $${\max_{\# Z = N} \frac{1}{N^2} \sum_{i,j = 1}^N d(z_i , z_j)  = \frac12},$$ and this maximum is achieved if and only if $Z$ is a  centrally symmetric set.
\item\label{M4} If $N$  is odd, $${\max_{\# Z = N} \frac{1}{N^2} \sum_{i,j = 1}^N d(z_i , z_j)  = \frac12 - \frac{1}{2N^2}},$$ and this maximum is achieved if and only if $Z$ can be represented as a union $ Z = Z_1 \cup Z_2$, where $Z_1$ is symmetric, while $Z_2$ lies on a two-dimensional hyperplane (i.e. on the same great circle) and satisfies $\displaystyle{\frac{1}{M^2} \sum_{z_i,z_j \in Z_2} d(z_i , z_j)  = \frac12 - \frac{1}{2M^2}}$  with  $M = \# Z_2$, i.e. $Z_2$ is a maximizer of the sum of geodesic distances on $\mathbb S^1$.
\end{enumerate}
\end{theorem}



Before we turn to the proof of the theorem, we briefly discuss the history of these questions. These results have been previously known in dimensions $d=1$ and $d=2$. Parts \eqref{M1} and \eqref{M2} in  $d=1$ were proved by Fejes-T\'oth \cite{FejesToth2} (and reproved later in \cite{d1} in relation to musical rhythms). Fejes-T\'oth also conjectured that the same holds for  for $d\ge 2$. In dimension $d= 2$ for even $N$ part \eqref{M3}of the theorem above was proved by Sperling \cite{Sperling}. Our Stolarsky principle \eqref{e.SPhem} makes the proof of this case very simple in all dimensions $d\ge 2$.   For odd $N$ Larcher \cite{Larcher} proved part \eqref{M4}  of Theorem \ref{t.geodpoints} in dimension $d=2$, however, we believe that there is a mistake in his proof (statement (ii) at the bottom of page  48).  We use a different geometric approach to prove  \eqref{M4}  in {\it{all dimensions}} $d\ge 2$, based  on part \eqref{M2} of Theorem \ref{t.geodpoints} and an interesting fact from combinatorial geometry known as the Sylvester--Gallai theorem.

 Thus, our   Stolarsky-type formula \eqref{e.SPhem} greatly simplifies the proof of these facts in dimensions one and two and, moreover, allows us to extend them to all dimensions $d\ge 2$. We now turn to the proof of Theorem \ref{t.geodpoints}.
 
 \begin{proof} Part \eqref{M1} of the theorem is now obvious since  the left-hand side of \eqref{e.SPhem} is non-negative and $\int\limits_{\mathbb S^{d}} \int\limits_{\mathbb S^{d}} d( x, y)  \, d\sigma (x)\, d\sigma (y) = \frac12$. 
 
 Part \eqref{M2} also follows easily from \eqref{e.SPhem}. Indeed, for every $x\in \mathbb S^d$ such that $x^\perp$ does not contain any points of $Z$ the minimal value of the integrand in the left-hand side of \eqref{e.SPhem}, i.e. the integrand in \eqref{e.dischem}, equals $0$ for even $N$ (if exactly half the points lie on either side of $x^\perp$), and is $\frac{1}{4N^2}$ for odd $N$ (if the numbers  of points on both sides of $x^\perp$ differ exactly by $1$). Obviously, configurations $Z$ for which this is achieved for each such $x \in \mathbb S^d$ are possible: e.g., $\lfloor N/2 \rfloor$  and $\lceil N/2 \rceil$ points in antipodal poles. Moreover, if  
 for some $x\in \mathbb S^d$ with $x^\perp \cap Z = \emptyset$ this condition is not satisfied, then it also fails on a small set of positive measure around $x$, hence $D_{L^2,\textup{hem}} (Z)$ is not minimal, and therefore  $\frac1{N^2} \sum d(z_i,z_j)$ is not maximized. \\
 
 To prove part \eqref{M3}, first observe that symmetric sets $Z$ trivially satisfy the condition of part \eqref{M2}. Now assume that for some $x\in \mathbb S^d$ the number of points of $Z$ located at $x$ and $-x$ is not the same. Consider a hyperplane passing through $x$, which contains no other points of $Z$. Perturbing it in opposite directions, we find that the difference of  number of points on either side changes by at least $2$, 
  i.e. cannot stay equal to  zero. Thus non-symmetric sets $Z$ with even number of points don't satisfy the condition of part \eqref{M2}, i.e. cannot maximize the sum of geodesic distances. \\
 
 We now turn to part \eqref{M4}. We shall rely on the Sylvester--Gallai theorem. In the Euclidean case it states the following: 
 {\it{if a finite set $Z$   in   $\mathbb R^d$ has the property that  for every two points of $Z$, the straight line passing through them contains at least one other point of $Z$, then all points of $Z$ lie on the same straight line.}}
 A spherical version of this theorem  also holds. 
  \begin{theorem}[Spherical Sylvester--Gallai]\label{t.sphSG}
 Assume that a set $Z$ of  $N$ points on the sphere $\mathbb S^d$ contains no antipodal points and satisfies the following condition: for every two points of $Z$, the great circle  passing through them contains at least one more point of $Z$. Then all points of $Z$ lie on the same great circle.
 \end{theorem} 
 For the history and several proofs of these theorems we refer the reader to the book \cite[pages 73 and 88]{PftB}.  Normally, these theorems are stated in dimension $d=2$, but higher dimensional extensions are simple. Indeed, for $Z \subset \mathbb S^d$, consider a copy of $\mathbb S^2$ which contains $z_1$, $z_2$, $z_3 \in Z$. The two-dimensional version of Theorem \ref{t.sphSG} applies, and hence $z_1$, $z_2$, $z_3$ lie on the same great circle. In the same manner, considering a copy of $\mathbb S^2$ containing this great circle and any other point $z_i \in Z$, we find that $z_i$ has to lie on the same great circle. 
 
 We are now ready to prove part \eqref{M4}. Assume that $N$ is odd. It follows from     \eqref{e.SPhem} and the proof of part \eqref{M2} that the maximal value of $\frac{1}{N^2} \sum d(z_i,z_j)$ is $\frac{1}{2} - \frac{1}{2N^2}$. Observe that adding a pair of antipodal points to $Z$ does not change maximality of $Z$, i.e.  $Z$ is a maximizer if and only if $Z\cup \{p, -p\}$ is a maximizer (with $N$ replaced by $N+2$). Indeed, since $d(x,p) + d(x,-p) =1$, it is easy to check that $$ \sum_{x,y \in Z \cup \{-p, p\} } d (x,y) =  \sum_{x,y \in Z  } d (x,y) + 2(N+1),$$ thus the second sum equals $\frac{N^2}{2} - \frac12$ if and only if the first sum is $\frac{(N+2)^2}{2} - \frac12$.  This immediately proves sufficiency of the condition in \eqref{M4}. Moreover, it shows that, in order to prove necessity, it is enough to  consider maximizers without antipodal points and to prove that they have to be contained in some great circle.
 
 Assume that  $Z$ with $\# Z = N$ maximizes  $\frac1{N^2} \sum d (z_i, z_j)$ and contains no pair of antipodal points. Consider two arbitrary points$z_1$, $z_2 \in Z$, and assume that no other point of $Z$ lies on the great circle defined by $z_1$ and $z_2$. Since $Z$ is finite, there exists a hyperplane containing $z_1$ and $z_2$, which does not contain any other points of $Z$. Since $z_1$ and $z_2$ are no antipodal, one can perturb the hyperplane in such a way that it does not touch other points of $Z$ and both points $z_1$ and $z_2$ end up on the same side of the hyperplane.  Perturbing in the opposite direction, we observe that the difference between the number of points on opposite side of the hyperplane changes by $4$, i.e. cannot stay equal to $\pm 1$, i.e., by part \eqref{M2},  $Z$ cannot be a maximizer. 
 
 We thus conclude that, for any two points of $Z$, at least one other point of $Z$ has to lie on the same great circle, i.e. the spherical Sylvester--Gallai theorem, Theorem \ref{t.sphSG}, applies. Hence $Z$ is contained in a great circle.
 \end{proof}


 \noindent{\emph{Remark:}} Observe that the one-dimensional maximizers of odd cardinality $N$, which arise in part \eqref{M4} of Theorem \ref{t.geodpoints}, are characterized by the condition that the sum of any $\lceil N/2 \rceil$ consecutive central angles defined by the points is at least $\pi$. In particular, any acute triangle is a maximizer  for $d=1$ and $N=3$.


Theorem \ref{t.geodpoints} demonstrates that the  situation is drastically different from the spherical  cap discrepancy  and  the sum of Euclidean distances. In the latter case, minimizing the $L^2$ spherical cap discrepancy (equivalently, maximizing the sum of Euclidean distances) leads to a rather uniform distribution of $Z$. 
In particular, for $d=1$ the sum is maximized by the vertices of a regular $N$-gon \cite{FejesToth}, and in higher dimensions maximizing distributions have to be well-separated \cite{St1}. The sum of {\emph{geodesic}} distances, however, may be maximized by very non-uniform sets, e.g. $N/2$ points in two antipodal poles. 

\section{Geodesic distance energy integral}\label{GeodInt}

The results of the previous section naturally  suggest 
a more general problem of finding  equilibrium distributions for the geodesic energy integral. Let $\mu$, $\nu$ be   Borel measures on $\mathbb S^d$. Define the geodesic distance energy integrals as
\begin{equation}
I (\mu, \nu) = \int\limits_{\mathbb S^d} \int\limits_{\mathbb S^d} d(x,y) \, d\mu (x) d\nu (y);\quad I(\mu) = I (\mu, \mu) =  \int\limits_{\mathbb S^d} \int\limits_{\mathbb S^d} d(x,y) \, d\mu (x) d\mu (y).
\end{equation}
Let $\mathcal M$ denote the set of non-negative Borel probability measures on $\mathbb S^d$ (i.e. $\mu (\mathbb S^d) =1$).  We shall be interested in the quantity
\begin{equation}
M = \sup_{\mu \in \mathcal M } I (\mu),
\end{equation}
as well as the maximizers of this expression, i.e.  the measures $\mu \in \mathcal M$ for which $I (\mu ) = M$ (the existence of maximizers follows easily from the weak$^*$-compactness of $\mathcal M$). The hemisphere Stolarsky principle, Theorem \ref{SPhem}, may be extended  to more general measures than the counting measure $\mu = \frac1{N} \sum_{i=1}^N \delta_{z_i}$.

\begin{theorem}[Hemisphere Stolarsky principle for general measures]\label{SPmu} Let $\mu$ be a Borel measure  on $\mathbb S^d$ with $\mu (\mathbb S^d) =1$. Then the following relation holds
\begin{equation}\label{e.SPmu}
\int\limits_{\mathbb S^d} \bigg( \mu \big( H(x) \big)  - \frac12 \bigg)^2 d\sigma (x)  = \frac12 \cdot \bigg( \frac12 - I (\mu) \bigg).
\end{equation}
\end{theorem}


\begin{proof} Notice that
\begin{equation}
\int\limits_{\mathbb S^d} \!\int\limits_{H(x)} d\mu (y) d\sigma (x) = \int\limits_{\mathbb S^d} \!\int\limits_{H(y)}  d\sigma (x)  d\mu (y) =  \frac12 \cdot \int\limits_{\mathbb S^d}   d\mu (y) = \frac12
\end{equation}
and, according to \eqref{2hem},
\begin{align}
\int\limits_{\mathbb S^d} \!\int\limits_{H(x)}  \!\int\limits_{H(x)} d\mu (y) d\mu (z)  d\sigma (x)  &  
=  \int\limits_{\mathbb S^d} \!\int\limits_{\mathbb S^d} \sigma \big( {H(y)\cap H(z) }  \big) \, d\mu (y)  d\mu (z) \\
&  = \int\limits_{\mathbb S^d} \!\int\limits_{\mathbb S^d} \frac12\cdot \big(1 - d(y,z) \big) \, d\mu (y)  d\mu (z)  = \frac12 - \frac12 I(\mu).
\end{align}
Using the two relations above we obtain
\begin{align}
\int\limits_{\mathbb S^d} \bigg( \mu \big( H(x) \big)  - \frac12 \bigg)^2 \, d\sigma (x)  &  = \int\limits_{\mathbb S^d} \bigg( \int\limits_{H(x)}  \!\int\limits_{H(x)} d\mu (y) d\mu (z)  -  \int\limits_{H(x)} d\mu (y)   + \frac14 \bigg)\,\, d\sigma(x)\\
& = \frac12- \frac12I(\mu)  - \frac12 + \frac14 = \frac12 \cdot \bigg( \frac12 - I (\mu) \bigg),
\end{align}
which proves \eqref{e.SPmu}. \end{proof}

Since the left-hand side of  identity \eqref{e.SPmu} is non-negative,
Theorem \ref{SPmu}  immediately yields a corollary about  the maximizers of $I(\mu)$:
\begin{corollary}\label{maximizers}
For any $\mu \in \mathcal M$, $I(\mu ) \le \frac12.$  Measures $\mu \in \mathcal M$, for which $I(\mu) = \frac12$, are exactly the measures which satisfy the following condition:
\begin{equation}\label{condition}
\mu \big( H(x) \big)  = \frac12 \,\,\, \textup{ for } \sigma\textup{-a.e. }\,\, x \in \mathbb S^d.
\end{equation}
\end{corollary}

It is very easy to see that  if the measure $\mu $ is symmetric, it is  a maximizer of $I(\mu)$, i.e. $I(\mu) = \frac12$. Indeed,  let $\mu^*$ be the reflection of $\mu$, i.e. $\mu^* (E) =  \mu (-E)$. It is easy to see that

\begin{align}
I (\mu, \mu^* ) & = \int\limits_{\mathbb S^d} \!\int\limits_{\mathbb S^d}  d(x,y) d\mu (x)  d\mu^* (y) =  \int\limits_{\mathbb S^d} \!\int\limits_{\mathbb S^d}  d(x, -y) d\mu (x)  d\mu (y) \\
& = \int\limits_{\mathbb S^d} \int\limits_{\mathbb S^d} \big( 1-  d(x, y)\big)  d\mu (x)  d\mu (y) = 1- I(\mu).
\end{align}
If, moreover $\mu$ is symmetric, i.e. $\mu^* = \mu$, then
\begin{equation}\label{e.symmax}
I(\mu) = I (\mu, \mu^*) = 1 - I(\mu) \,\,\, \Rightarrow \,\,\, I(\mu) = \frac12.
\end{equation}
Therefore, in particular, every symmetric measure $\mu \in \mathcal M$ satisfies \eqref{condition}.  The converse of this fact is  less obvious.



\begin{proposition}\label{l.hemsym} Assume that the measure $\mu \in \mathcal M$ satisfies the condition
\begin{equation}\label{condition1}
\mu \big( H(x) \big)  = \frac12 \,\,\, \textup{ for } \sigma\textup{-a.e. }\,\, x \in \mathbb S^d.
\end{equation}
Then the measure $\mu$ is symmetric, i.e. $\mu (E) = \mu (- E)$ for {every} Borel set $E \subset \mathbb S^d$.
\end{proposition}


We are unaware of an elementary proof of this seemingly simple statement. Our approach is based on spherical harmonics and Gegenbauer polynomials. We refer the reader to \cite{DX,Groemer} for background information. Let $w_\ld(t)=(1-t^2)^{\ld-\f12}$ with $\ld>0$.  Given $1\leq p<\infty$, we  denote by $L_{w_\ld}^p[-1,1]$ the space of all real  integrable functions $f$ on $[-1,1]$ with $\|f\|_{p,\ld}:=\Bl(\int_{-1}^1 |f(t)|^p w_\ld(t)\, dt\Br)^{1/p}<\infty$.
Every function $f\in L_{w_\ld}^1[-1,1]$ has  an expansion in terms of  Gegenbauer (ultraspherical) polynomials $C_n^\lambda$:
 \begin{equation}\label{1-1-16}f(t)\sim \sum_{n=0}^\infty \wh{f}(n,\ld) \f {n+\ld} \ld C_n^\ld(t),\  \  t\in [-1,1].\end{equation}
Let $\mathcal H_n (\sph)$ denote the space of spherical harmonics of degree $n$ on $\sph$, i.e. homogeneous harmonic polynomials of degree $n$ in $d+1$ variables.  We start with  an auxiliary lemma, which will also  be useful  in \S \ref{GenSt}.

 \begin{lemma}\label{FH} Let $\gamma $ be a signed Borel measure on $\mathbb S^d$ and $f\in L_{w_\ld}^2[-1,1]$ with $\lambda = \frac{d-1}{2}$.  Assume that
 \begin{equation}\label{cond}
 \int_{\sph} f(x\cdot y)\, d\gamma (y) = 0 \textup{ for } \sigma\textup{-almost every } x \in \mathbb S^d.
 \end{equation}
 Assume also that $\wh{f} (n,\lambda) \neq 0$ for some $n\ge 1 $. Then for every spherical harmonic of order $n$, $Y_n \in \mathcal H_n ( \mathbb S^d)$,
 \begin{equation}\label{orthog}
  \int_{\mathbb S^d} Y_n (y) d\gamma (y) = 0.
  \end{equation}
 \end{lemma}

 \begin{proof} The Funk-Hecke formula (see, e.g., Theorem 1.2.9 in \cite{DX}) states that
 \begin{equation}\label{e.FH}
  \int_{\mathbb S^d} f (x\cdot  y) Y_n (x)  d\sigma (x)  = \wh{f} (n,\lambda) Y_n (y).
  \end{equation}
 Therefore,
 \begin{align*}
  \int_{\mathbb S^d} Y_n (y) d\gamma (y) & =   \frac1{\wh{f} (n, \lambda)} \int_{\mathbb S^d}\int_{\mathbb S^d} f (x\cdot  y) Y_n (x) d\sigma (x) d\gamma (y)\\
   & =  \frac1{\wh{f} (n, \lambda)} \int_{\mathbb S^d} \bigg( \int_{\mathbb S^d} f (x\cdot  y)  d\gamma (y) \bigg) Y_n (x) d\sigma (x)  = 0.
 \end{align*}

 \end{proof}

\begin{proof}[Proof of Proposition \ref{l.hemsym}]
Let $\mu^*$ be a reflection of $\mu$, defined by $\mu^* (E)  = \mu (-E)$, and set $\gamma  = \mu - \mu^*$.  Condition \eqref{condition1} then implies that $$\gamma \big( H(x) \big) =  \int_{\sph} {\bf 1}_{(0,1]}(x\cdot y)\, d\gamma (y) = 0 \,\,\, \textup{ for } \sigma\textup{-a.e. }\,\, x \in \mathbb S^d.
 $$
 The Gegenbauer coefficients $\widehat{f} (n,\lambda)$ of  the function $f(t) = {\bf 1}_{(0,1]} (t)$  are non-zero for odd $n$ (Lemma 3.4.6 in \cite{Groemer}). Therefore, according to Lemma \ref{FH}, relation \eqref{orthog} holds for all odd $n$. But for even values of $n$ it obviously holds, since in this case  $Y_n \in \mathcal H_n (\sph)$ is an even function, and $\gamma$ is antisymmetric.
  Therefore, $\int_{\mathbb S^d} f (y) d\gamma (y) = 0$ for every polynomial $f$, and hence for each $f \in C (\sph)$, which implies that $\gamma =0$. Hence $\mu = \mu^*$, i.e. $\mu$ is symmetric.
  \end{proof}

From the above discussion we obtain the following characterization of the maximizers of $I(\mu)$:
\begin{theorem}\label{1max}
For a measure $\mu \in \mathcal M$, $$ I(\mu ) = \sup_{\gamma \in \mathcal M} I (\gamma) = \frac12$$ if and only if $\mu$ is centrally symmetric.
\end{theorem}

This behavior of $I(\mu )$  goes in sharp contrast with the behavior of the seemingly similar energy integral $\int\limits_{\mathbb S^d} \int\limits_{\mathbb S^d} \| x- y\|  \, d\mu (x) d\mu (y)$. It is  known  \cite{Bjorck} that the unique maximizer of this energy integral is $\mu = \sigma$, the uniform distribution on $\mathbb S^d$. In this sense the behavior of $I(\mu)$ is more similar  (albeit still different) to that of $\int\limits_{\mathbb S^d} \int\limits_{\mathbb S^d} \| x- y\|^2  \, d\mu (x) d\mu (y)$ which is maximized by any measure with center of mass is at the origin, which may  be easily seen from the relation
\begin{equation}\label{e.norm2}
\int\limits_{\mathbb S^d} \int\limits_{\mathbb S^d} \| x- y\|^2  \, d\mu (x) d\mu (y) = 2 - 2\cdot  \bigg\| \int\limits_{\mathbb S^d} x d\mu (x) \bigg\|^2 .
\end{equation}
It is thus natural to analyze  energy integrals with  general powers $\delta >0$.

\subsection{Geodesic distance energy integrals with exponent $\delta  >0$.}
We would like to understand which measures $\mu \in \mathcal M$ maximize the energy
\begin{equation}
I (\mu) = \int\limits_{\mathbb S^d} \int\limits_{\mathbb S^d} \big( d(x,y) \big)^\delta \, d\mu (x) d\mu (y),
\end{equation}
for $\delta >0$, and how  maximizers depend on $\delta$.

While the geodesic distance energy integral is a novel object, such  integrals with Euclidean distances are well investigated.  For an extensive study of the energy integrals $\int\limits_{F} \int\limits_{F} \| x- y\|^\delta  \, d\mu (x) d\mu (y)$, $\delta >0$,  see \cite{Bjorck}. Specialized to the case $F = \mathbb S^d$, these results are formulated below.
\begin{theorem}[Bjorck, \cite{Bjorck}]\label{Bjo}
For  $\delta>0$, define the energy integral
\begin{equation}
I_E (\mu) = \int\limits_{\mathbb S^d} \int\limits_{\mathbb S^d} \| x- y\|^\delta  \, d\mu (x) d\mu (y)\red{.}
\end{equation}
The maximizers of this energy integral over  $\mu\in \mathcal M$ (Borel probability measures on $\mathbb S^d$) can be characterized as follows:
\begin{enumerate}[(i)]
\item\label{E1} $0 < \delta < 2$: the unique maximizer of $I_E (\mu)$ is $\mu = \sigma$ (the normalized surface measure).
\item\label{E2} $\delta = 2$:  $I_E (\mu)$ is maximized if and only if the center of mass of $\mu$ is at the origin.
\item\label{E3}  $\delta > 2$: $I_E (\mu)$ is maximized if and only if  $\mu = \frac12 (\delta_p + \delta_{-p})$, i.e. the mass is equally  concentrated at two antipodal poles.
\end{enumerate}
\end{theorem}

The proof of part \eqref{E1} uses potential analysis, in particular, the semigroup property of the Riesz potentials; part \eqref{E2} is explained in \eqref{e.norm2}; and part \eqref{E3} is almost trivial.

We  observe that there is a ``breaking point" $\delta = 2$ in the behavior of maximizers of the Euclidean  energy integral. Surprisingly,  for the seemingly similar geodesic distance integral this critical value  is different: $\delta =1$. We have the following theorem:
\begin{theorem}\label{GeodMax}
For  $\delta>0$ let $I (\mu)$ be the geodesic distance energy integral
\begin{equation}
I (\mu) = \int\limits_{\mathbb S^d} \int\limits_{\mathbb S^d} \big( d(x,y) \big)^\delta  \, d\mu (x) d\mu (y)
\end{equation}
The maximizers of this energy integral over  $\mu\in \mathcal M$  can be characterized as follows:
\begin{enumerate}[(i)]
\item\label{G1} $0 < \delta < 1$: the unique maximizer of $I (\mu)$ is $\mu = \sigma$ (the normalized surface measure).
\item\label{G2} $\delta = 1$:  $I (\mu)$ is maximized if and only if  $\mu$ is centrally symmetric .
\item\label{G3} $\delta > 1$: $I (\mu)$ is maximized if and only if  $\mu = \frac12 (\delta_p + \delta_{-p})$, i.e. the mass is equally  concentrated at two antipodal poles.
\end{enumerate}
\end{theorem}

Part \eqref{G1} is proved in the companion paper of the authors \cite{BD2} through extensive  analysis of spherical harmonics expansions. Part \eqref{G2} is the result of Theorem \ref{1max} above, which is a consequence of  the hemisphere Stolarsky principle \eqref{e.SPmu} and  is contained in  Lemma \ref{l.hemsym}.  The proof of part \eqref{G3} is quite simple: since $d(x,y) \le 1$, we have for $\delta >1$
\begin{align*}
I (\mu ) \le   \int\limits_{\mathbb S^d} \int\limits_{\mathbb S^d}  d(x,y)   \, d\mu (x) d\mu (y) \le \frac12.
\end{align*}
The first inequality turns into an equality when $\mu \times \mu \big\{  (x,y): \, d(x,y) = 0 \textup{ or } 1 \big\} =1$, while the second bound becomes exact when $\mu$ is symmetric, according to part \eqref{G2}. This readily implies that $\mu =  \frac12 (\delta_p + \delta_{-p})$.

This peculiar effect (that geodesic distance energy behaves differently from its Euclidean counterpart) has been noticed in dimension $d=1$, i.e. on the circle, in \cite{BHS}, where the  one-dimensional case of parts \eqref{G1} and \eqref{G3} of the above theorem have been proved. In \cite{BD2}, the follow-up to the present paper, we conduct a more detailed analysis of the geodesic distance energy (including negative powers and logarithmic energy).

\subsection{Average case integration error} Finally, the right-hand side of the hemisphere Stolarsky principle \eqref{e.SPhem} also yields the average-case squared integration error on $C(\mathbb S^d)$ with respect to the law of the {\emph{hemisphere Gaussian process}} introduced in \cite{BL}. This is a mean-zero  Gaussian process  $G$ on $\mathbb S^d$, which is defined by $\mathbb E G^2_x  = \frac14$, $\mathbb E \big( G_x - G_y \big)^2 = \| {\bf 1}_{H(x)}- {\bf 1}_{H(x)} \|_2^2 = d(x,y)$, i.e.  its covariance is given by $\mathbb E G_x G_y = \frac14 - \frac12 d(x,y)$.  It induces a Gaussian measure on the space of continuous functions $C (\mathbb S^d)$, which we also denote by $G$. Then the average-case integration error with respect to $G$  is equal to the $L^2$ hemisphere discrepancy.
\begin{theorem}\label{l.average}
Let $Z = \{ z_1, z_2, \dots, z_N \} \subset \mathbb S^d$. The following holds
\begin{align}
\int\limits_{C(\mathbb S^d)}  \bigg|  \frac1{N} \sum_{i=1}^N f(z_i)  - & \int\limits_{\mathbb S^d} f(x)d\sigma(x)  \bigg| ^2  \, dG(f)  \\    & = \frac12   \cdot \left(\, \frac12 - \frac{1}{N^2} \sum_{i,j = 1}^N d(z_i , z_j) \right)  =  [D_{L^2,\textup{hem}} (Z)]^2 .
\end{align}
\end{theorem}

\begin{proof} We note that
\begin{equation}
\int\limits_{C(\mathbb S^d)} f(x) f(y) dG(f) = \mathbb E G_x G_y = \frac14 -\frac12 d(x,y),\,\,
\end{equation}
and thus
\begin{equation}
\int\limits_{C(\mathbb S^d)} \bigg( \int\limits_{\mathbb S^d}  f(x) d\sigma (x) \bigg)^2 \, dG(f)  = \int\limits_{\mathbb S^d} \int\limits_{\mathbb S^d} \bigg(  \frac14 -\frac12 d(x,y)\bigg) d\sigma (x) d\sigma (y) = 0,
\end{equation}
i.e. $ \int\limits_{\mathbb S^d}  f(x) d\sigma (x) = 0$ for $G$-a.e. $f\in C (\mathbb S^d)$. Therefore, for the average case  integration error we obtain
\begin{align}
\int\limits_{C(\mathbb S^d)}  \bigg|  \frac1{N} \sum_{i=1}^N f(z_i)  - & \int\limits_{\mathbb S^d} f(x)d\sigma(x)  \bigg| ^2  \, dG(f)  = \int\limits_{C(\mathbb S^d)}  \bigg|  \frac1{N} \sum_{i=1}^N f(z_i)    \bigg| ^2  \, dG(f)\\
& = \frac1{N^2} \sum_{i,j =1}^N \bigg( \frac14 - \frac12 d(z_i,z_j) \bigg)  = \frac12 \cdot \bigg( \frac12 -  \frac1{N^2} \sum_{i,j =1}^N d(z_i,z_j) \bigg),
\end{align}
which is exactly the right-hand side of the hemisphere Stolarsky principle \eqref{e.SPhem}. \end{proof}

We remark that the first result of this type has been obtained in \cite{WandCo} for the anchored $L^2$ discrepancy on $[0,1]^d$  and the average-case integration error with respect to the Wiener sheet measure. \\

\section{Positive definite functions and generalized Stolarsky principle}\label{GenSt}

In this section we take a more general look on energy minimization and Stolarsky principle. For any Borel measure $\mu$ on $\sph$ and a bounded or non-negative Borel measurable function $F$ on $[-1,1]$, we define the energy integral     $$ I_F (\mu ):=\int_{\sph}\int_{\sph} F (x\cdot y)\, d\mu(x) \, d\mu(y).$$ As before, let $ \mathcal B$ denote the class of {finite}  Borel {signed}  measures on $\sph$.

We start with a simple observation which shows that, while $I_F (\mu)$ is quadratic in $\mu$, 
it behaves linearly near $\sigma$.

\begin{lemma}\label{musigma}
For any  $F$ be a bounded or non-negative Borel measurable function on $[-1,1]$, and a signed measure  $\mu \in \mathcal B$ with $\mu (\sph)=1$, the following relation holds:
\begin{equation}\label{e.musigma}
I_ F(\mu)  - I_F (\sigma)  = I_F ({\mu - \sigma} ).
\end{equation}
\end{lemma}

\begin{proof}
The proof followed from  a simple observation that, due to  rotational invariance, for any $x\in \sph$ $$ \int_{\sph} F (x\cdot y) d\sigma (y) = \int_{\sph} \int_{\sph} F (x\cdot y) d\sigma (x) d\sigma (y) = I_F (\sigma ),$$ i.e., the left-hand side is independent of $x\in \sph$. Therefore,
\begin{align*}
I_ F({\mu-\sigma})  & = I_F  (\mu ) + I_F (\sigma ) - 2 \int_{\sph} \int_{\sph} F (x\cdot y) d\sigma (x) d\mu (y)\\
& = I_F (\mu )  + I_F (\sigma ) - 2 I_F (\sigma)  \int_{\sph} d\mu (y) = I_F (\mu) - I_F (\sigma ).
\end{align*}
\end{proof}




Next, we recall the concept of positive definite functions on the sphere. A function $F \in C[-1,1]$ is called {\emph{positive definite}} on the sphere $\sph$ if for any set of points $Z = \{ z_1, ..., z_N \} \subset \sph$, the matrix $\big[ F (z_i\cdot z_j )\big]_{i,j=1}^N$ is positive semidefinite, i.e.
\begin{equation}\label{defpd}
\sum_{i,j} F(z_i \cdot z_j) c_i c_j \ge 0
\end{equation}
 for all $c_i \in \mathbb R$. We denote the class of positive definite functions by $\Phi_d$. This class admits several different characterizations.

\begin{proposition}\label{tfae}
For a function  $F \in C[-1,1]$ the following conditions are equivalent:
\begin{enumerate}[(i)]
\item\label{i} $F$ is positive definite on $\mathbb S^d$, i.e. $F\in \Phi_d$.
\item\label{ii} For $\lambda =\frac{d-1}{2}$, all Gegenbauer coefficients of $F$ are non-negative, i.e.
\begin{equation}
\widehat{F} (n, \lambda) \ge 0\,\, \textup{ for all } n\ge 0.
\end{equation}
\item\label{iii} For any signed measure $\mu \in \mathcal B$ the energy integral  is non-negative: $I_ F (\mu)\ge 0$.
\item\label{iv} There exists a function $f \in L^2_{w_\lambda}[-1,1]$ such that
\begin{equation}\label{3-1}
    F(x\cdot y)=\int_{\sph} f(x\cdot z) f(z\cdot y)\, d\s(z),\   \ x, y\in \sph,
\end{equation}
i.e. $F$ is the spherical convolution of $f$ with itself.
\end{enumerate}
\end{proposition}

We shall briefly outline the proof of this proposition. The equivalence of \eqref{i} and \eqref{ii} is a celebrated theorem of Schoenberg \cite{Schoen}. In addition it is known that the Gegenbauer expansion of $F\in \Phi_d$ is absolutely summable. Since \eqref{defpd} states that $I_F \big(\sum c_i \delta_{z_i} \big) \ge 0$, obviously \eqref{iii} implies \eqref{i}. The converse implication is proved by a standard argument based on the compactness of $\sph$ and the weak-$*$ density of the linear span of Dirac masses in $\mathcal B$. Finally, the equivalence of \eqref{ii} and \eqref{iv} can be established by defining $f$ through the identity $\big( \widehat{f} (n, \lambda) \big)^2  = \widehat{F} (n, \lambda)$. 
Absolute summability of the Gegenbauer series of $F \in \Phi_d$   will guarantee that $f \in L^2_{w_\lambda} [-1,1]$. For more details on positive definite functions, see Chapter 14 in \cite{DX}.\\



Condition \eqref{iii} above suggests that the property of being positive definite is related to energy minimization. We show that this is indeed the case. We shall need to make a technical assumption that $I_\sigma (F) \ge 0$ (in view of  \eqref{iii}, it is necessary for $F\in \Phi_d$). But, since adding a constant to $F$ does not effect minimizing energy over $\mathcal M$, this assumption is easily removable. We first prove the following theorem.

\begin{theorem}\label{sigmamin}
Assume that $F \in C[-1,1]$ and $I_\sigma (F) \ge 0$. Then $\sigma$ is a minimizer of $I_F (\mu)$ over $\mathcal M$ (probability Borel measures on $\mathbb S^d$) if and only if $F \in \Phi_d$.
\end{theorem}

\begin{proof}
The sufficiency follows easily from \eqref{e.musigma} and condition \eqref{iii} of Proposition \ref{tfae}. Indeed, if $F\in \Phi_d$, then for any $\mu \in \mathcal M$
\begin{equation}
I_F (\mu) - I_F (\sigma) = I_F (\mu - \sigma) \ge 0, \,\, \textup{ i.e. } \,\,  I_F (\mu) \ge I_F (\sigma).
\end{equation}

We now prove the necessity. Assume that $\sigma$ is a minimizer of $I_F (\mu)$ over $\mu \in \mathcal M$. We first state an auxiliary lemma:
\begin{lemma}\label{zero}
Let  $F \in C[-1,1]$ with $I_\sigma (F) \ge 0$. Assume that $\sigma$ is a minimizer of $I_F (\mu)$ over $\mu \in \mathcal M$. Then  for any measure $\gamma \in \mathcal B$ with total mass zero, $I_F (\gamma) \ge 0$.
\end{lemma}
Assuming the lemma, it is easy to finish the proof of the theorem. Indeed, let $\mu \in \mathcal B$ be an arbitrary signed measure with total mass one. Then according to \eqref{e.musigma} and the lemma above $$ I_F (\mu ) - I_F (\sigma) = I_F (\mu - \sigma) \ge 0,$$ since $(\mu - \sigma ) (\mathbb S^d) =0$. Therefore, $I_F (\mu ) \ge I_F (\sigma) \ge 0$, and by part \eqref{iii} of Proposition \ref{tfae}, $F$ is positive definite.
\end{proof}

\begin{remark} We would like to observe that along the way we have proved that if $\sigma$ is ia minimizer of $I_F$ over $\mathcal M $, the set of positive measures of mass one, it is also a minimizer over the class of all {\emph{signed}} measures of total mass one.  This is not necessarily the case in other settings. In particular, for the integral over the ball $$ \int_{\mathbb B^{d+1}} \int_{\mathbb B^{d+1}} \| x-y \| d\mu (x) d \mu (y),$$ according to \cite{Bjorck} the unique maximizer over $\mathcal M$ is $\sigma$, while in the case of signed measures the maximizer does not exist \cite{Hin}.
\end{remark}

It remains to prove Lemma \ref{zero}
\begin{proof}[Proof of Lemma \ref{zero}]
Assume, on the contrary, that for some $\gamma \in \mathcal B$ with $\gamma (\sph) = 0$, we have $I_ F (\gamma) < 0 $. We shall smooth $\gamma $ out by considering, for $\varepsilon >0$,  the function $\rho_\varepsilon (x) = \frac{\gamma \big( C (x,1-\varepsilon) \big)}{\sigma \big( C (x,1-\varepsilon) \big)}$ and defining $d\gamma_\varepsilon (x) = \rho_\varepsilon (x) d\sigma (x)$,where as before $C(x,t)=\{ z\in \mathbb S^d:\, x\cdot z > t\} $ is the spherical cap. The measure $\gamma_\varepsilon$ has total mass zero, since, letting $m_\varepsilon = \sigma \big( C (x,1-\varepsilon) \big)$, we have  $$ \int_{\sph} d\gamma_\varepsilon (x) = \frac{1}{m_\varepsilon}  \int_{\sph} \int_{C(x,1-\varepsilon) }  {d\gamma (y)}  d\sigma (x)  =   \frac{1}{m_\varepsilon}\int_{\sph} \bigg( \int_{C(y,1-\varepsilon)} d\sigma (x) \bigg)  {d\gamma (y)}     = 0.$$

Next, we claim that, for $\varepsilon $ small enough, $I_F (\gamma_\varepsilon) < 0$. This follows from
\begin{align*}
{I_F} (\gamma_\varepsilon) & = \int\limits_{\sph} \int\limits_{\sph} \frac{1}{m_\varepsilon^2} \int\limits_{C(x,1-\varepsilon) } \int\limits_{C(y,1-\varepsilon) }   F(x\cdot y)   {d\gamma (u)}    {d\gamma (v)} d\sigma (x)   d\sigma (y)\\
& = \int\limits_{\sph} \int\limits_{\sph} \left[ \frac{1}{m_\varepsilon^2} \int\limits_{C(u,1-\varepsilon) } \int\limits_{C(v,1-\varepsilon) }   F(x\cdot y)  d\sigma (x)   d\sigma (y) \right]  {d\gamma (u)}    {d\gamma (v)} \longrightarrow I_F (\gamma)< 0
\end{align*}
as $\varepsilon  \rightarrow 0$, since the expression inside the brackets converges to $F(u \cdot v)$ uniformly in $u$ and $v$. This proves the claim.

It is also  easy to see that the density $\rho_\varepsilon$ is a bounded function. Therefore, there exists a constant $c>0$ such that the measure $\mu = c\gamma_\varepsilon + \sigma$, i.e. $d\mu (x) = (1 + c\rho_\varepsilon (x)) d\sigma (x)$, is non-negative. Hence $\mu \in \mathcal M$. Since $\sigma$ minimizes $I_F$ over $\mathcal M$, by \eqref{e.musigma} we have $$ 0 \le  I_F (\mu ) - I_F (\sigma) = I_F (\mu - \sigma ) = c^2 I_F (\gamma_\varepsilon),$$ which contradicts the fact that $I_F (\gamma_\varepsilon) < 0$.
\end{proof}

Since $I_F (\sigma ) = \widehat{F} (0, \lambda)$ for $\lambda = \frac{d-1}2$, we can easily remove the assumption $I_F (\sigma) \ge 0$ in Theorem \ref{sigmamin}.

\begin{corollary}
Assume that $F \in C[-1,1]$. Then $\sigma$ is a minimizer of $I_F (\mu)$ over $\mathcal M$  if and only if $F + C \in \Phi_d$ for some constant $C\in \mathbb R$ or, equivalently, if $\widehat{F} (n, \lambda) \ge 0$ for all $n \ge 1$.
\end{corollary}

\vskip3mm

 We now turn to the generalization of Stolarsky principle for positive definite functions. Assume that $F\in\Phi_d$, $\lambda = \frac{d-1}2$, and the function  $f\in L_{w_\ld}^2[-1,1]$ is such that \eqref{3-1} is satisfied.

For a  non-negative Borel probability measure $\mu $ on $\mathbb S^d$ we  define the $L^2$ discrepancy of $\mu$ with respect to $f$ as
\begin{align*}
    D_{L^2, f}(\mu) =\left(\int_{\sph}\Bl|\int_{\sph} f(x\cdot y)\, d\mu(y) - \int_{\sph} f(x\cdot y)\, d\s(y) \Br|^2\, d\s(x)\right)^{\f12}.
\end{align*}
The $L^2$ discrepancy of a finite point-set $Z \subset \sph$ is simply
\begin{equation}
 D_{L^2, f}(Z) =  D_{L^2, f}\Big( \frac1N {\sum_{j=1}^N} \delta_{z_j}\Big) =  \left(\int_{\sph}\Bl| \frac{1}{N} \sum_{i=1}^N f (z_i)  - \int_{\sph} f(x\cdot y)\, d\s(y) \Br|^2\, d\s(x)\right)^{\f12}.
\end{equation}
Notice that various choices of $f$ recover different geometric notions of discrepancy, although this object is more general. We now prove a general version of the Stolarsky principle, which connects the energies with respect to $F$ to the $L^2$ discrepancy built upon $f$.


\begin{theorem}[Generalized Stolarsky principle]\label{thm-3-2}  Let $\mu \in \mathcal B $ be a signed Borel probability measure on $\mathbb S^d$ with total mass $\mu (\sph) =1$ and let $F \in \Phi_d$ with $f$ as in \eqref{3-1}. Then
\begin{align}\label{stolmu}
    I_F (\mu  )  - I_F (\sigma  ) =  D_{L^2, f}^2(\mu) .
\end{align}
In particular,  in the case of $\mu = \frac1{N} \sum_{i=1}^N \delta_{z_i}$, this relation becomes
\begin{align}\label{stolmu1}
    \frac{1}{N^2} \sum_{i,j=1}^N F (z_i \cdot z_j)  - \int\limits_{\sph} \int\limits_{\sph} F( x\cdot y) d\sigma (x) d\sigma (y)   =  D_{L^2, f}^2(Z) .
\end{align}

 \end{theorem}

\begin{proof} According to the definition of $ D_{L^2, f}(\mu)$,  \eqref{e.musigma}, and \eqref{3-1}, we have
\begin{align*}
  D^2_{L^2, f}(\mu) & = \int_{\sph}\Bl( \int_{\sph} f(x\cdot y)\, d(\mu - \sigma) (y)  \Br)^2\, d\s(x) \\
  & =  \int_{\sph} \int_{\sph} \int_{\sph} f(x\cdot y) f(x\cdot z)  d(\mu - \sigma) (y)  d(\mu - \sigma) (z) d\s (x) \\
  & =    \int_{\sph} \int_{\sph} F(y \cdot z)   d(\mu - \sigma) (y)  d(\mu - \sigma) (z) = I_F ({\mu - \s}  ) =   I_F (\mu )  - I_F (\sigma  ) .
\end{align*}

\end{proof}

This approach brings up several novel points. First of all, in most contexts Stolarsky identity arises from the notion of the $L^2$ discrepancy, which in turn dictates the specific form of the interaction potential $F$. Theorem \ref{thm-3-2}, on the other hand, allows one to go in the opposite direction: starting with the potential $F\in \Phi_d$, one can produce a natural notion of discrepancy, for which the Stolarsky principle holds. The precise form of the function $f$,  defined through the identity $\big( \widehat{f} (n, \lambda) \big)^2  = \widehat{F} (n, \lambda)$,  cannot be made explicit in most cases (in fact, many different choices of $f$ corresponding to the same $F$ can be constructed by changing the signs of the coefficients $\widehat{f} (n, \lambda) $). However, this does not prevent one from being able to obtain estimates for $D_{L^2,f} (\mu)$. In \cite{BD2} (Theorem 4.2, part (ii)) we prove that
\begin{align}
   C_d \min_{1\leq k\leq c_d N^{1/d}} \wh{F}(k,\ld)\leq   \inf_{\# Z = N } {D}_{L^2, f}^2 (Z) \leq  N^{-1}\max_ {0\leq \theta \leq c_d' N^{-\f1d}} \bl( F(1)-   F(\cos \theta)\br).\label{3-2}\end{align}
Hence, e.g., lower bounds can be proved using information  about either $F$ or $f$. In \cite{BD2} we use these  estimates to give an alternative proof of the spherical cap discrepancy bounds \eqref{Beck},  and employ \eqref{stolmu1} {to} obtain sharp asymptotic behavior of the difference between discrete and continuous energies, $E_F (Z) - I_F (\sigma)$, as the number of points $N\rightarrow \infty$ both in the case of Riesz energy, $F(x\cdot y ) = \| x-y \|^{\delta}$  (recovering results of \cite{Wagner,SK,Brauchart}), and the geodesic distance energies, $F (x\cdot y)  = \big( d(x,y) \big)^\delta$, introduced in this paper.

Here we concentrate on the applications of the Stolarsky principle \eqref{stolmu} to characterizing minimizers of $I_F$.  Since $D^2_{L^2, f}(\mu) \ge 0$,  identity \eqref{stolmu} gives yet another proof that for  $F \in \Phi_d$, the uniform measure  $\sigma$ is a minimizer of $I_F (\mu   )$ over $\mathcal M$ (in fact, over all {\emph{signed}} Borel probability measure on $\mathbb S^d$ with total mass $\mu (\sph) =1$). Furthermore, the generalized Stolarsky identity \eqref{stolmu} also allows one to characterize those $F\in \Phi_d$ for which $\sigma$ is the {\it{unique}} minimizer of $I_\mu (F)$.

\begin{theorem}\label{sunique} Let $F \in {C[-1,1]}$. Then
 $\sigma$ is the unique minimizer of $I_F (\mu )$ if and only if $ \wh{F} (n, \lambda) >0$ for each $n\ge 1$.
\end{theorem}





 We shall need a lemma  which is a  simple corollary of Lemma \ref{FH}   and the density of polynomials in $C (\mathbb S^d)$ -- compare it to the proof of Proposition \ref{l.hemsym}.

 \begin{lemma}\label{corol} Let $\gamma$ be a signed Borel measure on  $\mathbb S^d$ with $\gamma (\mathbb S^d) = 0$ and  let $f\in L_{w_\ld}^2[-1,1]$.  Assume that  condition \eqref{cond} of Lemma \ref{FH} is satisfied, i.e.  \begin{equation}
 \int_{\sph} f(x\cdot y)\, d\gamma (y) = 0 \textup{ for } \sigma\textup{-almost every } x \in \mathbb S^d,
 \end{equation} and $\wh{f} (n,\lambda) \neq 0$  for all $n \ge 1$.  Then $ \gamma = 0$.
 \end{lemma}

We are now ready to prove Theorem \ref{sunique}:

 \begin{proof} {We start with the proof of sufficiency.  Without loss of generality, we may assume that $\wh{F}(0; \ld)>0$.} Assume that for each $n\ge 1$, we have $ \wh{f} (n; \lambda) = ({\wh{F} (n; \lambda)})^{1/2} \neq 0$. Let $\mu$ be a minimizer of $I_\mu (F)$, i.e. $I_\mu (F) = I_\sigma (F)$. Therefore,  the Stolarsky principle \eqref{stolmu} implies that $  D^2_{L^2, f}(\mu) = 0$, i.e. $\int_{\sph} f(x\cdot y)\, d(\s-\mu) (y) = 0$ for $\sigma$-almost every $x$. Then by  Corollary \ref{corol}, $\sigma - \mu =0$. Hence $\sigma$ is the unique minimizer of $I_\mu (F)$. \\

Conversely, assume that  $\wh{F} (n_0, \lambda)   \leq  0$ for some $n_0\ge 1$. Let $Y_{n_0}$ be a spherical harmonic of degree $n_0$ with $\|Y_{n_0}\|_2=1$.
Then, for $\varepsilon >0$ small enough, the measure $ d \mu = \big(1 + \varepsilon Y_{n_0}(x) \big) d\sigma  \in \mathcal M$, and by the  Funk-Hecke formula \eqref{e.FH}, we have
\begin{align*}
    I_\mu (F)&=\int_{\sph}\int_{\sph} F(x\cdot y) \Bl(1+\va Y_{n_0}(x)\Br)\Bl(1+\va Y_{n_0}(y)\Br)\, d\s(x) \, d\s(y)\\
    &=I_\sa(F)+\va^2 \int_{\sph}\int_{\sph} F(x\cdot y) Y_{n_0}(x) Y_{n_0}(y)\, \, d\s(x) d\s(y)\\
    &+2\va \int_{\sph}\int_{\sph} F(x\cdot y) Y_{n_0}(x)\, d\s(x)\, d\s(y)\\
    &=I_\sa(F) +\va^2 \wh{F}(n_0;\ld)\leq I_\sa(F),
\end{align*}
with equality being valid only if $\wh{F}(n_0;\ld)=0$. This is impossible since  $d\s$ is the unique minimizer of $I_\mu$(F).

 \end{proof}


The relations between positive definite functions and energy minimization on the sphere are well known \cite{Schoen,SK}. Here we have attempted to give an essentially self-contained exposition with minimal references to ultraspherical expansions and a special emphasis on the novel role of Stolarsky principle. It has come to our attention that  similar ideas have been explored also in \cite{DH}.

\section{Stolarsky principle for spherical wedges and slices.}\label{WedSli}

The spherical wedge as defined in \cite{BLd} is the subregion of  $\mathbb S^d$ between two hyperplanes: for $x$, $y \in \mathbb S^d$
\begin{equation}
W_{xy} = \{  z \in \mathbb S^d: \, \textup{sign} ( x \cdot z ) \neq  \textup{sign} ( y \cdot z ) \} ,
\end{equation}
in other words it is the collection of all points $z \in \mathbb S^d$ such that the hyperplane $z^\perp$ separates $x$ and $y$. It is  easy to see that 
(compare to  \eqref{2hem})
\begin{equation}
\sigma (W_{xy} )  = {\mathbb \sigma} \big\{ \textup{sign} (x\cdot z) \neq  \textup{sign} (y\cdot z ) \big\} =  \sigma \big(H (x) \Delta H(y)\big)=d(x,y).
\end{equation}
(this is also a simple instance of the Crofton formula in integral geometry).
For a finite set of vectors $Z=\{ z_1, z_2, ..., z_N\}$ on the sphere $\mathbb S^{d}$, we define the Hamming distance between the  points $x$, $y\in \mathbb S^{d}$ as
\begin{equation}
d_H (x,y) : = \frac1{N} \cdot {\# \big\{ z_k \in Z:\, \textup{sign} (x\cdot z_k) \neq  \textup{sign} (y\cdot z_k ) } \big\} = \frac{\# \{ Z \cap W_{xy}  \} }{N},
\end{equation}
i.e.   the proportion of those hyperplanes $z_k^\perp$ that {\it{separate}} the points $x$ and $y$.  Therefore, the quantity
\begin{equation}
\Delta_Z (x,y) := d_H (x,y)  -    d(x,y) = \frac1{N} \sum_{i=1}^N {\bf 1}_{W_{xy}} (z_i)  - \sigma (W_{xy})
\end{equation}
is precisely the discrepancy of $Z$ with respect to $W_{xy}$. This quantity arises  naturally in {\emph{one-bit compressed sensing}}, uniform tessellations of the sphere, as well as dimension reduction and almost isometric embedding  results (e.g., one-bit analogs of the Johnson--Lindenstrauss lemma), see \cite{BL,BLd,PV} for more details.

We define the  $L^2$ wedge discrepancy
\begin{equation}
D_{L^2,\textup{wedge} }  (Z) =  \big\| \Delta_Z (x, y) \big\|_2 
 =
  \left( \,\, \int\limits_{\mathbb S^d}\int\limits_{\mathbb S^d} \bigg( \frac1{N} \sum_{i=1}^N {\bf 1}_{W_{xy}} (z_i)  - \sigma (W_{xy}) \bigg)^2 d\sigma(x) d\sigma(y)  \right)^{\frac12}.
\end{equation}

The analog of the Stolarsky principle for wedges has been proved by the first author and M. Lacey in \cite{BLd}:
\begin{theorem}[Stolarsky principle for wedges]\label{SPwedge}
\begin{equation}\label{e.SPwedge}
[D_{L^2,\textup{wedge} }  (Z) ]^2
    =   \frac{1}{N^2} \sum_{i,j = 1}^N \bigg( \frac12 - d(z_i , z_j) \bigg)^2   -  \int\limits_{\mathbb S^{d}} \int\limits_{\mathbb S^{d}} \bigg( \frac12 - d( x, y) \bigg)^2  \, d\sigma (x)\, d\sigma (y)  .
\end{equation}
\end{theorem}

This theorem implies that, in order to minimize $D_{L^2,\textup{wedge} }  (Z)$ one should minimize the discrete energy with the potential $\big( \frac12 - d(x,y)\big)^2$, i.e. make the vectors $z_k$ as orthogonal as possible on the average. We would like to point out  the strong similarity between this discrete energy and another similar quantity, the {\emph{frame potential}} introduced in \cite{BF}:
\begin{equation}
FP (Z)  = {\frac1{N^2}} \sum_{i,j =1}^N  ( z_i \cdot z_j )^2 .
\end{equation}
A  finite set $Z = \{ z_1, \dots, z_N \} \subset \mathbb S^d$ is called a {\emph{tight frame}} if and only if there exists a constant $C>0$ such that for any vector $x\in \mathbb R^{d+1}$
\begin{equation}
\| x \|^2  = C \sum_{i=1}^N (x \cdot z_j )^2 .
\end{equation}  It was proved in \cite{BF} that $Z$ is a tight frame if and only if $Z$ is a minimizer of the frame potential $FP (Z)$. \\

Similarly, we  define the spherical slices: for $x$, $y \in \mathbb S^d$
\begin{equation}
S_{xy} = \{ z \in \mathbb S^d :\, x\cdot z > 0,\, y \cdot z < 0 \},
\end{equation}
i.e. a slice is half of a wedge. The $L^2$ slice discrepancy naturally is
\begin{equation}\label{e.dslice}
D_{L^2,\textup{slice} }  (Z)  =     \left( \,\, \int\limits_{\mathbb S^d}\int\limits_{\mathbb S^d} \bigg( \frac1{N} \sum_{i=1}^N {\bf 1}_{S_{xy}} (z_i)  - \sigma (S_{xy}) \bigg)^2 d\sigma(x) d\sigma(y)  \right)^{\frac12}.
\end{equation}
This discrepancy has been previously considered in \cite{Blu}, however the Stolarsky principle in this setting is new:
\begin{theorem}[Stolarsky principle for slices]\label{SPslice}
\begin{equation}\label{e.SPslice}
4 [D_{L^2,\textup{slice} }  (Z) ]^2   =       \frac{1}{N^2} \sum_{i,j = 1}^N \big( 1 - d(z_i , z_j) \big)^2   -  \int\limits_{\mathbb S^{d}} \int\limits_{\mathbb S^{d}} \big( 1- d( x, y) \big)^2  \, d\sigma (x)\, d\sigma (y)    .
\end{equation}
\end{theorem}
\begin{proof} This proof is very similar to the proof of  \eqref{e.SPwedge}, see \cite[Theorem 1.21]{BLd}.  Recall that $ S_{xy} = \{ z \in \mathbb S^d :\, x\cdot z > 0,\, y \cdot z < 0 \}$ and $\sigma (S_{xy} ) = \frac12 d(x,y)$. Using  the definition \eqref{e.dslice} we obtain

\begin{align}
\nonumber  & [D_{L^2,\textup{slice} }  (Z) ]^2  =
  \frac{1}{N^2} \sum_{i,j = 1}^N \int\limits_{\mathbb S^d}\! \int\limits_{\mathbb S^d}     {\bf 1}_{S_{xy}} (z_i) \cdot  {\bf 1}_{S_{xy}} ({z_j}) \, d\sigma (x) \, d\sigma (y)\\
 & \qquad  - \frac{2}{N}  \sum_{k =1}^N \int\limits_{\mathbb S^d} \!\int\limits_{\mathbb S^d}   {\bf 1}_{S_{xy} } (z_k)  \cdot \sigma (S_{xy}) \,\, d\sigma (x) \, d\sigma (y)  +  \int\limits_{\mathbb S^{d}} \!\int\limits_{\mathbb S^{d}}  \sigma (S_{xy})^2 \, d\sigma (x) \, d\sigma (y).
\end{align}\vskip2mm

 It is easy to see that $z_i$, $z_j \in S_{xy}$ if and only if $x\in S_{z_i, -z_j}$ and $y \in S_{-z_i, z_j}$. Since
\begin{equation}
\sigma \big( S_{  \pm z_i, \mp z_j} \big) = \frac12 d ( \pm z_i, \mp z_j ) = \frac12 \big(1 - d (z_i, z_j ) \big) ,
\end{equation}
we find that
\begin{equation}
\int\limits_{\mathbb S^d} \int\limits_{\mathbb S^d}   {\bf 1}_{S_{xy}} ({z_j}) \cdot  {\bf 1}_{S_{xy}} (z_i)  d\sigma (x) d\sigma (y) = \sigma  \big( S_{z_i,-z_j} \big) \cdot  \sigma  \big( S_{-z_i,z_j} \big) = \frac14 \cdot \big( 1- d(z_i,z_j) \big)^2.
\end{equation}

Notice that by rotational invariance the double integral in the second term does not depend on the choice of $z_k \in \mathbb S^d$ and therefore it can be replaced by the average over $z\in \mathbb S^d$:
\begin{align}
\int\limits_{\mathbb S^d} \!\int\limits_{\mathbb S^d}   &{\bf 1}_{S_{xy} } (z_k)  \cdot \sigma (S_{xy})  \,\, d\sigma (x) \, d\sigma (y)  = \int\limits_{\mathbb S^d} \!\int\limits_{\mathbb S^d} \!\int\limits_{\mathbb S^d}   {\bf 1}_{S_{xy} } (z)  \cdot \sigma (S_{xy}) \,\, d\sigma (x) \, d\sigma (y) \, d\sigma (z) \\
&  =   \int\limits_{\mathbb S^d} \!\int\limits_{\mathbb S^d} \bigg[ \!\int\limits_{\mathbb S^d}   {\bf 1}_{S_{xy} } (z) \, d\sigma (z) \bigg]\,  \sigma (S_{xy}) \,\, d\sigma (x) \, d\sigma (y)  = \int\limits_{\mathbb S^d} \!\int\limits_{\mathbb S^d}   \sigma (S_{xy})^2  \,\, d\sigma (x) \, d\sigma (y).
\end{align}
Since $\sigma (S_{xy} ) = \frac12 d(x,y)$, it follows that
\begin{align}
[D_{L^2,\textup{slice} }  (Z) ]^2  & =    \frac{1}{4N^2} \sum_{i,j = 1}^N \big( 1- d(z_i,z_j) \big)^2  -  \frac14 \int\limits_{\mathbb S^d} \!\int\limits_{\mathbb S^d}   \big( d(x,y) \big)^2  \,\, d\sigma (x) \, d\sigma (y)\\
& =    \frac{1}{4N^2} \sum_{i,j = 1}^N \big( 1- d(z_i,z_j) \big)^2  -  \frac14 \int\limits_{\mathbb S^d} \!\int\limits_{\mathbb S^d}   \big( 1 - d(x,y) \big)^2  \,\, d\sigma (x) \, d\sigma (y),
\end{align}
which proves Theorem \ref{SPslice}. \quad \end{proof}

The  value of the integral $\int\limits_{\sph} \int\limits_{\sph}  \big( d(x,y) \big)^2  \,\, d\sigma (x) \, d\sigma (y) $,  
  which arises in this theorem,  will be computed in the next section.

We also  note   that the optimal order of magnitude both for $D_{L^2,\textup{slice} }  (Z) $  \cite{Blu} and $D_{L^2,\textup{wedge} }  (Z) $ \cite{BLd} satisfy the same bounds as   the spherical cap discrepancy \eqref{Beck}.


\section{Appendix: mean-square geodesic distance}\label{App}
The following integral arises in the formulations of Stolarsky principles for wedges \eqref{e.SPwedge} and slices \eqref{e.SPslice}:
\begin{equation}
V_d  = 	\int\limits_{\mathbb S^d} \!\int\limits_{\mathbb S^d}   \big( d(x,y) \big)^2  \,\, d\sigma (x) \, d\sigma (y),
\end{equation}
hence we compute it and  examine its properties. A standard calculation yields that
\begin{equation}
V_d = \frac{1}{\pi^2} \cdot  \frac{\omega_{d-1}}{\omega_d} \int_{0}^\pi \phi^2  \sin^{d-1} \phi \, d\phi.
\end{equation}
Applying the recursive relation \cite{GR}   
\begin{align}
\int x^m \sin^n x \, dx & = \frac{x^{m-1} \sin^{n-1} x}{n^2 } \, \big(m \sin x  - nx \cos x \big) + \\
& + \frac{n-1}{n} \int x^m \sin^{n-2} x \,dx - \frac{m(m-1)}{n^2} \int x^{m-2} \sin^n x \, dx
\end{align}
with $m=2$ and $n=d-1$, as well as the facts that
\begin{equation}
\frac{\omega_{d-1}}{\omega_d} = \frac{d-1}{d-2}  \cdot \frac{\omega_{d-3}}{\omega_{d-2}}
\end{equation}
and
\begin{equation}
\int_0^\pi \sin^{d-1} \phi \, d\phi = \frac{\sqrt{\pi} \Gamma (d/2)}{\Gamma \big( (d+1)/2\big)},
\end{equation}
 we obtain the recursive relation
\begin{equation}
V_d  = V_{d-2}  - \frac{2}{\pi^2 (d-1)^2}.
\end{equation}
Together with  simple identities $V_1 = \frac{1}{3}$ and $V_2 = \frac{1}{2} - \frac{2}{\pi^2}$  (or even $V_0 = \frac12$) this yields
\begin{lemma}
For odd values of $d\ge 1$
\begin{equation}
V_d =  \frac13 - \frac{2}{\pi^2} \sum_{k=1}^{{(d-1)}/{2}} \frac{1}{(2k)^2 },
\end{equation}
while for even values of $d\ge 2$
\begin{equation}
V_d = \frac12 -  \frac{2}{\pi^2} \sum_{k=1}^{{d}/{2}} \frac{1}{(2k-1)^2 }.
\end{equation}
\end{lemma}
Since $\displaystyle{\sum_{k=1}^\infty \frac1{(2k)^2}  = \frac{\pi^2}{24}}$ and $\displaystyle{\sum_{k=1}^\infty \frac1{(2k-1)^2}  = \frac{\pi^2}{8}}$, we find that $$\displaystyle{\lim_{d\rightarrow \infty} V_d = \frac14},$$ which is consistent with the concentration of measure phenomenon (``most points" on the high-dimensional sphere are nearly orthogonal).

Notice that this confirms the result of Theorem \ref{GeodMax} that, unless $d=0$, the uniform distribution $\sigma$ is {\it{not}} a maximizer of $I (\mu)  = \int\limits_{\mathbb S^d} \!\int\limits_{\mathbb S^d} \big( d(x,y) \big)^2  \,\, d\mu (x) \, d\mu (y)$, since for $\mu = \frac12 {\delta_p }+  \frac12 {\delta_{-p}}$ we have $I (\mu) = \frac12$.

\end{document}